\documentclass[12pt]{article}
\usepackage{amsmath,amsthm,amssymb}
\usepackage[final]{graphicx}
\usepackage[pdftex,pdfpagelayout=OneColumn]{hyperref}

\newtheorem{thm}{\sc Theorem}

\newtheorem{lem}[thm]{\sc Lemma}
\newtheorem*{smallworldtheorem}{\sc Small World Theorem}
\newtheorem*{cosmologicaltheorem}{\sc Cosmological Theorem}
\newtheorem*{growththeorem}{\sc Growth Theorem}
\newtheorem*{frequencytheorem}{\sc Frequency Theorem}
\newtheorem*{splittingtheorem}{\sc Splitting Theorem}

\theoremstyle{definition}
\newtheorem*{defn}{\sc Definition}

\theoremstyle{remark}
\newtheorem*{rem}{\sc Remark}

\newtheoremstyle{algorithm}	
     {\topsep}  
     {\topsep}  
     {}         
     {}         
     {\bfseries}
     {}         
     {\newline} 
     {}         
\theoremstyle{algorithm}
\newtheorem*{alg}{\sc Algorithm}

\newcommand{\V}{\vskip0.4in}

\begin{document}
\begin{center}
\Large
Look, There's More to Say about Conway's Look and Say Sequence
\end{center}

\begin{flushright}
Gregory P.~Dresden \\
Washington and Lee University \\
Lexington, VA 24450\\
\verb+dresdeng@wlu.edu+

\vspace{2 mm}

Jacob Siehler \\
Washington and Lee University\\
Lexington, VA 24450\\
\verb+siehlerj@wlu.edu+
\end{flushright}



\vspace{10 mm}

Consider the many different ways of describing the string
\[
S_0=1, 2, 1, 3, 5, 5.
\]
We might say that this string has ``two 1's, one 2, one 3, and two
5's", and we could write down this description as a new string:
\[
S_1= 2, 1, \,1, 2, \,1, 3, \,2, 5.
\]
One naturally thinks about continuing this process, producing a
sequence of strings $S_2, S_3, S_4$, and so on, each string describing
the previous one.  This sequence is often called a {\bf Counting
Sequence}, and as seen in \cite{BF, SS} 
all such sequences are ultimately periodic.  (The
above sequence becomes constant at the string $S_3=3, 1, 1, 2, 3,
3, 1, 5$).
 
Another way to describe $S_0$ is to indicate the number of times a
digit appears without writing the digit itself.  So, our description
of $S_0$ now becomes
\[
T_1= 0, 2, 1, 1, 0, 2,
\]
indicating zero 0's, two 1's, one 2, one 3, zero 4's, and two 5's.
(It's essential to include those 0's in $T_1$ as placeholders.) The
following table gives a few iterations of this method:

\vskip0.2in

\begin{tabular}{r|cccccc} 
  &  0's   &   1's   &   2's  &    3's  &    4's   &   5's \\
\hline 
{$T_1$}  &    0  &       2 &       1   &     1  &       0  &     2  \\
{$T_2$}   &   2   &      2  &      2    &    0   &      0   &    0 \\
{$T_3$}    &  3    &     0   &     3     &   0    &     0    &   0 
\end{tabular}

\vskip0.2in

As one might expect, all such sequences of strings  eventually
become periodic, as seen in \cite{K}, \cite{MW}, and \cite{SE}.
(The above sequence becomes periodic, with period 2, at the strings
$T_7= 3, 1, 1, 1, 0, 0$ and $T_8= 2, 3, 0, 1, 0, 0$.) In the
literature, these are called {\bf Self-Descriptive} or {\bf
Co-Descriptive} sequences.  See also \cite {PS}, \cite{G} for
variations on this theme involving infinite strings.  The problem
of finding fixed points under this iterative method has also appeared
in such disparate sources as Marilyn vos Savant's infamous ``Ask
Marilyn" column in Parade Magazine \cite{VS}, and the in-flight
magazine for the now-defunct airline TWA \cite{MW}.

In this paper, we consider a third way of describing the string
$S_0$. The following method
was introduced by John Conway in \cite{C} but it dates
back to at least 1977 and a meeting between the Dutch and British
teams at the International Mathematical Olympiad \cite{MPW},
and possibly has an even older history.  Sometimes called the 
{\bf Look and Say} or {\bf Counting Consecutive Appearances} method,
it describes $S_0$ as ``one 1, one 2, one 1, one 3 and two 5's",
which we write as
\[
U_1= 1, 1, 1, 2, 1, 1, 1, 3, 2, 5
\]
If we drop the commas and concatenate the digits,
we get the following sequence of strings:
\begin{eqnarray*}
U_1 &=& 1112111325 \\
U_2 &=& 311231131215 \\
U_3 &=& 13211213211311121115
\end{eqnarray*}
Some computer scientists (and aficionados of computing history) will recognize
this type of sequence as a {\bf Morris sequence}, named after
the NSA cryptographer Robert Morris and mentioned as a puzzle in a popular
book from 1989 \cite{Stoll}.

We note that Conway's process is simply an iteration of the {\bf Run-Length
encoding} method of data compression familiar to
computer scientists. 
A {\bf run} is a (sub)string of just one repeated digit, such as
$11111111$, and run-length encoding would map that to $91$ (indicating a run
of nine $1$'s). 
Likewise, the more generic string {\tt xxxxxxyyyyyyyyyyyy}
with two runs would be mapped to $6${\tt x}$12${\tt y}.

Conway demonstrated in \cite{C} that when repeating this Look and Say process, 
almost any initial
string
eventually breaks apart (or splits), in the following sense.
For $S$ a string, we say that $S$ {\bf splits} if we can write 
$S = A.B$ (the period simply indicates
concatenation)
in such a way that the later iterates $S_{n}$ of $S$ 
all have the form
$A_{n}.B_{n}$, where $A_n$ is the $n$th iteration of $A$
under this Look and Say method,
and likewise $B_n$ for $B$.
If a string does not split, we say that it is {\bf
irreducible}. In the example given above, the strings $U_1, U_2, U_3$ can
be rewritten with the following splits (indicated by a '.'):
\begin{eqnarray*}
U_1 &=& 1112.11132.5 \\
U_2 &=& 3112.311312.15 \\
U_3 &=& 132112.1321131112.1115
\end{eqnarray*}

At this point it's a good idea to take a pencil and work out the
next few iterations of the look-say process on $U_3$, to see how the
$2$'s in this string form a natural `barrier\rq\  which justifies the
splitting points we have indicated. You may also begin to see why the
three strings in the splitting of $U_1$ are irreducible.  
We won't prove that here, but for reassurance you
can check one of the references that follow.

Conway discovered that there were 92 special and irreducible strings
(or {\bf elements}) such as $1112$ and $11132$ that kept appearing and
reappearing under this iterative process. In a fit of whimsy, Conway 
named these 92 elements after the 92 chemical elements of the periodic
table, and in his ``Cosmological Theorem" he stated that any initial
string (except for the fixed string 22 and the empty string
$\epsilon$) will eventually split into all 92 elements (each repeated
many times) along with, possibly, a few ``exotic" elements.
Furthermore, the length of the string under this iterative process
grows as $\lambda^n$, for $\lambda = 1.303577\dots$ the single
positive real root of a 71 degree polynomial.

John Conway and Richard Parker wrote down a proof of the Cosmological
Theorem which they then lost; Mike Guy found a shorter proof which
he then lost as well.  A complete (and computer-assisted) proof was
later published by D. Zeilberger \cite{EZ} and improved by R. A.
Litherland \cite{Li}.  For a review of Conway's results and an
extension to the topic of infinite self-referential strings, see O.
Mart\'{i}n's recent paper \cite{M} on Look-and-Say Biochemistry.


We call attention to the following subtle point. 
Conway's Look and Say method would read this
string
\[
V_0= 5 ,5, 5, 5, 5, 5, 5, 5, 5, 5
\]
as ``ten 5's", and would write this description as
\[
V_1= 10,5
\]
which becomes 
\begin{eqnarray*}
 V_2 & =& 1, 10, 1, 5 \qquad (\mbox{instead of } 1, 1, 1, 0, 1, 5)\\
V_3 &=& 1, 1, 1, 10, 1, 1, 1, 5
\end{eqnarray*}
It is this feature (obscured by the concatenation
of digits in many of the cited papers) that leads us to
say that Conway's method operates in ``base $\infty$".
Let's see what would happen instead if we worked
in base $10$, concatenating the digits.
Then, we would have
\begin{eqnarray*}
V_1^* &=& 105 \\
V_2^* &=& 111015 \qquad (\mbox{instead of } V_2 = 11015)\\
V_3^* &=& 31101115, 
\end{eqnarray*}
which appears, at first glance, to differ from the earlier sequence
$V_0, V_1, V_2, V_3, \dots$, yet (surprisingly enough) it ends
up behaving exactly as predicted by Conway's theory: after a finite
number of ``steps\rq\rq\ (that is, iterations), 
all 92 of Conway's common elements will appear in
every string, and the growth rate will be
$\lambda$. 

So, does the choice of base matter at all? It turns out that for
any base greater than base $3$, all strings will decay exactly as
predicted by Conway's theorem, with the same $92$ common elements
and the same $\lambda$. This leaves only bases $2$ and $3$ to
consider, and for these two cases the behavior is indeed a bit
different than in Conway's description.

The binary case was studied by T. Silke in 1993 \cite{Si}. He found
the binary world to be a much smaller world with only eleven common
elements, not all of which appear in every decayed string. He also
found that while strings still grow exponentially in length, the
rate of growth is different than in Conway's theory; the base $2$
rate of growth is $\lambda_2 = 1.4655\dots$, the single real root
of $x^3-x^2-1$.

Base 3 is the only case left to consider, and we are pleased to
report that it has much to interest us. For example, in base $3$
we find three fixed strings: $11110$, $11112$, and Conway's $22$.
There are 21 other common strings akin to Conway's chemical elements,
giving 24 in all, and as expected every starting string eventually
decays into combinations of those 24.  What is particularly interesting
(and unique to base 3) is that some starting strings eventually
produce all 24 common strings among their descendents, while other
starting strings are only able to produce a special subset of the
24. From this point forward, we will be working exclusively in base 3. 


\subsection*{The strange world of base-3 subatomic particles.}

As mentioned above, Conway named his 92 common strings (also known
as {\bf elements}) after the
chemical elements (hydrogen, helium, etc) of the periodic table.
Following Conway's example, yet keeping in mind that ours is a
smaller (base 3) world, we will call our 24 common strings {\bf particles}, 
and we will name them after
the elementary particles (electrons, muons, up quarks, etc)
of subatomic physics. These particles divide naturally into the 
three groups called
{\bf fermions}, {\bf bosons}, and {\bf neutrinos}, as
seen here in {\sc Table} \ref{base3}:

\begin{table}[!hbt]
\small
\begin{tabular}{cc@{\ \ \ \ \ \ \ \ }cc@{\ \ \ \ \ \ \ \ }cc}
\multicolumn{2}{c}{\bf Fermions} & \multicolumn{2}{c}{\bf Bosons\ \ \ \ \ \ \ } & \multicolumn{2}{c}{\bf Neutrinos} \\
\hline\hline
 & & 211 & (Ph) & 22 & (Ne) \\
 & & 1221 & (Gl) &  11110 & (Nm)   \\
 & & 112211 & (Wb) & 11112 & (Nt) \\
 & & 12221 & (Zb) &  \\
 & & 2 & (H) &  \\
  10 & (E)  & 	12 & (Se) &  & \\
1110 & (M)  & 	1112 & (Sm) &  & \\
110 & (U)  & 	112 & (Su) &  & \\
2110 & (D)  & 	2112 & (Sd) &  & \\
122110 & (S)  & 122112 & (Ss) &  & \\
11222110 & (C)& 11222112 & (Sc) &  & \\
22110 & (B)  & 	22112 & (Sb) &  & \\
222110 & (T)  & 222112 & (St) & &
 \end{tabular}
\caption{The 24 base{-}3 particles}\label{base3}
\end{table}

\V

We note that this table does not contain every irreducible
string! Instead, in keeping with Conway's definitions,
we are only listing the {\bf common} particles, those
which show up in sufficiently decayed strings. Consider,
for example, the obviously irreducible string 1. Its decay chain is
as follows:
\begin{eqnarray*}
&1&\\
&11&\\
&21&\\
&1211&\\
&111221&\\
&10.12211&\\
&1110.112221&\\
&10.110.2110.221& = \mbox{E.U.D.Ph}
\end{eqnarray*}
We see that the string 1 (and also 11 and 21 and so on) eventually
decays into four of our common particles, E (for Electron), U (for
Up quark), D (for Down quark), and Ph (for Photon).  On the other
hand, the string 1 itself will never show up as one of the irreducible
components in a sufficiently old string, which is why it is not
included in our table of 24 particles. For a further explanation of our particular
(physics-themed) 
nomenclature for these particles, see the Appendix. 

We note also that a bit of thought will indicate that the $0$'s in the 
strings above
will serve as a `barrier\rq, just as with the $2$'s in the strings
$U_1, U_2, U_3$ discussed earlier, and thus justifying our splitting.
We give complete details on how to split a base-$3$ string in our {\sc Splitting
Theorem} to follow, but for now, this lemma will suffice. 

\begin{lem}\label{zero}
A string can always be split at the point where $0$ is followed by
a non-$0$ digit.
\end{lem}


\subsection*{Main Results.}


Very briefly, the results that follow can be summed up as ``Eventually,
everything is almost entirely made of fermions.''  To make more
precise statements possible, we need a few definitions.  First, say
that a string $S$ is {\bf eventually common} if some iterate $S_n$
splits as a concatenation of some or all of the 24 common, named
particles introduced in the previous section.  We also find the
following definition convenient:

\begin{defn}
An {\bf essential ancient string} is a string that contains (i) 16
digits or fewer; (ii) no runs of a single digit of length 4 or more;
and (iii) at most a single 0, and that, if it occurs, in the final
position of the string.
\end{defn}

We consider these strings ``essential'' because of the following
theorem. The reason for calling them ``ancient'' will emerge in due
course.

\begin{smallworldtheorem}
If every essential ancient string is eventually common, then every
string is eventually common.
\end{smallworldtheorem}

We call this the ``Small World'' theorem because it gives us a
finite (and, in fact, tractable) search space which determines
the long-term behavior of every string in the ternary universe.
With computer assistance, then, we will also provide a proof of the following:

\begin{cosmologicaltheorem}
Every essential ancient string decays into a collection of common
particles in at most 10 iterations. Consequently, every string is
eventually common.
\end{cosmologicaltheorem}

Finally, a study of the decay pattern for common particles will allow us 
(by the end of this paper) to
deduce the following two results on long-term behavior:

\begin{growththeorem}
For a string $A$ which does not split as a concatenation of neutrinos,
the string length of the iterates $A_n$ grows at an exponential
rate:
\[\lim_{n\to\infty}\hbox{Length}(A_n)/\hbox{Length}(A_{n-1})=\lambda,\]
where $\lambda\approx 1.3247\ldots $ is the real root of the
polynomial $x^3-x-1$.
\end{growththeorem}

\begin{frequencytheorem}
For a string $A$ which does not split as a concatenation of neutrinos,
for sufficiently large $n$, the iterates $A_n$ will contain only a
vanishing proportion of bosons and neutrinos.  Furthermore, the
relative frequency of fermions in the iterates will tend to the
following percentages:

\begin{tabular}{l|l|l|l|l|l|l|l}
 E& M& D& B& U& S& T& C \\
 \hline
 18.50\% & 13.97\% & 13.97\% & 13.97\% & 10.54\% & 10.54\% & 10.54\% & 7.96\% \\
\end{tabular}

\end{frequencytheorem}

\subsection*{The Small World Theorem}

We begin with a few lemmas on the length of single-digit runs in a
string.

\begin{lem}\label{eli}
Suppose all runs in $A$ are of finite length.  Then, after a finite
number of iterations (depending on the size of $A$), all descendents of
$A$ will have runs of length seven or less.
\end{lem}

We note that the ``finite number of iterations" mentioned above is very
finite; in practice, three steps (that is, three iterations) will suffice.  Only for strings
with runs of length more than $10^{230}$ would four or more steps
be required.

\begin{proof}
Let $m>7$ be the length of the longest run in $A$.  We will show
that all runs in the first iterate $A_1$ are of length strictly
less than $m$ (so, after a certain number of steps, we can get down
to runs of length seven or less).  Write $A$ as:
\[
A =  \quad    .....x\ \ \underbrace{yyyyyy...y}_{\mbox{\footnotesize $n_1$ times}}
        \ \  \underbrace{zzzzzz...z\phantom{y}\!\!\!}_{\mbox{\footnotesize $n_2$ times}}....
\]

with $n_1 , n_2 \leq  m$. Let's now look at $A_1$:
\[
A_1 =  \quad    .....x\ (n_1\mbox{ in base 3})\ y\ 
                       (n_2\mbox{ in base 3})\ z...
\]
The highest number of $y$'s will occur when ``$n_1$ in base
3" and ``$n_2$ in base 3" are written all in $y$'s.  How many $y$'s
will we get?  The number of digits in the base-3 representation of
$n_1$ is bounded by $1+\log_3{n_1}$, and likewise for $n_2$, and
so adding in one for the ``$y$" itself we get a run of no more than
$(3+\log_3{n_1}+\log_3{n_2})$ $y$'s, which itself is bounded by
$3+2\log_3{m}$.  So long as $m>7$, it's easy to show that $3+2\log_3{m}
< m$, which means all runs in $A_1$ are of shorter length, as desired.
\end{proof}

\begin{lem}\label{elii}
Given a string $A$ with all runs of length  $ \leq 7$, then $A_2$
(and all its descendents) has at most one consecutive $0$, at most
three consecutive $2$'s, and at most four consecutive $1$'s.
\end{lem}

\begin{proof}
$A_1$ can't have 000, as this can only be parsed as either x00 0's
(implying at least nine 0's in $A$) or x000 y's (implying at least
twenty-seven y's in $A$).  Likewise, $A_2$ can't have 00, again by
parsing.  $A_1$ (and thus $A_2$)  can't have 2222, as this parses
as either 222 2's, or 22 2's, or maybe 2 2's and 22 y's, each
implying a run of at least eight in $A$.

Finally $A_2$ can't have 11111, as the most parsimonious parsing
possible would be 11 1's and 11 y's, implying at least four (non-1)
y's in $A_1$, which by our work above is impossible.
\end{proof}

Let's define a string to be {\bf mature} if both the string and its
parent satisfy the conclusions of Lemma \ref{elii} on the length
of runs of digits.  In the language of Lemma \ref{elii}, the
descendent $A_3$ will be mature since both $A_3$ and its parent
$A_2$ satisfy the statement of the lemma.

At this point, we need to discuss the behavior of the mu and tau
neutrinos (11110 and 11112) in mature strings.  If $A$ is sufficiently
old, then the only possible parsing of these two neutrinos is as
11 1's and 1 0, and 11 1's and 1 2, respectively.  Thus, these both
can only arise as children of themselves, and can not be created
through the decay of other particles or strings.

The mu neutrino ending in $0$, since it's inside a mature string, will 
not be followed by another $0$ and so (by {\sc Lemma} \ref{zero}) will
always split the string.  The tau neutrino is more subtle in its behavior.
Consider the decay of the (non-mature) string 11112 11112 1:
\begin{eqnarray*}
11112 & 11112\ 1  & \\
11112 & 11112\ 11  & \\
11112 & 11112\ 21  & \\
11112 & 1112211 &  \\
11112 & 10\ 12221 &  =\mbox{Nt.E.Zb}
\end{eqnarray*}
We see that the left 11112 does indeed remain a subatomic particle,
while the right 11112 eventually decays into an electron (E) and Z
boson (Zb), thanks to the disrupting presence of that additional 1
on the far right. This motivates the following lemma.

\begin{lem}\label{eliii}
Let $A$ be an irreducible, non-empty, mature string.  Then, after
at most  three iterations, $A$ decays into a mixture of common particles
and particles with all runs of length  $\leq 3$.
\end{lem}

\begin{rem}
We call a string {\bf ancient} if the string is mature and has no
runs of length four or more.  From our discussion of the mu and tau
neutrinos, we know that once a string becomes ancient, it remains
ancient.  By the way, if you've been keeping track, you'll note
that all reasonable strings (with runs of length $\leq 10^{230}$)
become ancient in at most nine steps.
\end{rem}

\begin{proof}[Proof of Lemma \ref{eliii}]
$A$ is mature, so the only possible run of length 4 is 1111.  If
$A$ ends in 1111, then $A_1$ ends in exactly three 1's.  If $A$
contains 11110, then since $A$ is irreducible it must equal 11110
and thus is common.  So, suppose $A$ contains 11112.  If $A$ ends
in 11112, then $A$ either \emph{is} 11112, or the final 11112 can
be split off.  So, suppose $A$ contains .....11112J..... where $J$
is some non-empty string.  If $J$ begins with a 2, then $A_1$ has
only three 1's (and thus becomes and stays ancient).  So, we are
left with the case where $J$ doesn't begin with 2, doesn't split,
and is inside a mature string.  There are only a few choices for
the initial part of $J$ (12, 112, 1112, 1112, 1110, etc.) and it
is not difficult to try out each of them in turn and discover that
either $J_1$ or $J_2$ begins with a 2, meaning that by the third
step the run of four 1's in $A$ has disappeared.
\end{proof}


\begin{proof}[Proof of Small World Theorem]
First, we show that every ancient string decays in four steps into
strings of length  $\leq 16$.  


Let $A$ be an irreducible, ancient string of length $>16$.  We can
write A in runs as follows, where (1's) represents a run of one,
two, or three 1's (a run of four 1's isn't possible in an ancient
string), and likewise for the symbol (2's).

\begin{center}\qquad A = (possibly some 1's)(2's)(1's)$\cdots$(2's)(possibly some 1's)(possible 0)
\qquad \qquad (*)\end{center}

Each step changes the (2's) into either 12 or 22 or 102, and the
(1's) become either 11 or 21 or 101.  Thus, in a given step the
number of distinct runs of 2 increases by at most one (occuring
only if A begins with 112..., which produces 21... the next step)
Our plan is to show that by the fourth step, the string has split
into substrings, each with no more than two runs of 2's (and thus
each being of size $\leq 16$).

Now, since A has more than 16 digits and is of the form  (*) , it
must have a subsequence of the following form:

\begin{tabular}[t]{rcccccc}
... & G & H & I & J & K & ... \\
... & (2's) & (1's)  &  (2's) &    (1's) &   (2's) & ...
\end{tabular}

If I is a run of three 2's, then on the first step this string splits
when I becomes 10.2 (the ``." is added to emphasize the splitting
of I).  The same holds for G and K so we assume these all have only
one or two 2's.

If I is a run of two 2's, there are three cases to consider, depending
on the length of J and K.  First, the two cases where J is a run
of one 1:

\begin{tabular}[t]{r|rrrr}
 & H & I & J & K \\
 \hline\hline
step 0 & (1's) & 22 & 1 & 2 \\
\hline
step 1 & (1's) & 22 & 11 & 12 \\
step 2 & (1's) & 22 & 10.1 & \\ 
\end{tabular}
	\ \ \ \ \ \ \ \ \ \ 
\begin{tabular}[t]{r|rrrr}
 & H & I & J & K \\
 \hline\hline
step 0 & (1's) & 22 & 1 & 22 \\
\hline
step 1 & (1's) & 22 & 11 & 22 \\
step 2 & (1's) & 22 & 21 &     \\
step 3 & (1's) & 10.2 &  &    
\end{tabular}

Second, the one case where J is a run of two 1's:

\begin{tabular}{r|rrrr}
 & H & I & J & K \\
 \hline\hline
step 0 & (1's) & 22 & 11 & (2's) \\
\hline
step 1 & (1's) & 22 & 21 & \\
step 2 & (1's) & 10.2 & & 
\end{tabular}

If I is a run of a single 2, then I must be flanked by 11 on both
the left and the right, because a mature substring of form 2121 or
1212 doesn't parse.  So, our substring has the form:

\begin{tabular}[t]{ccccccc}
... & G     & H  & I & J  & K     & ... \\
... & (2's) & 11 & 2 & 11 & (2's) & ...
\end{tabular}

There are various possibilities for the behavior of J and its
descendents; here is the worst possible case:

\begin{tabular}[t]{r|rrrrr}
 & G & H & I & J & K \\
 \hline\hline
step 0  &   (2's)  & 11 & 2    & 11 & 22 \\
\hline 
step 1  &   (2's)  & 1 1 & 2 2 & 1  & 22  \\
step 2  &   (2's)  & 1  & 22   & 11 & 22 \\
step 3  &   (2's)  & 11 & 22 2  & 1  & \\
step 4  &          &    & 10.2 &    & \\
\end{tabular}

So, by step 4, our initial string A has split (perhaps in several
places) into smaller strings.  We must consider the annoying
possibility that additional runs of 2's were generated in A or its
substrings over those four steps.  This could happen in two locations:
at the beginning of A, or at the beginning of one of its splittings.
Now, a new 2 is born at the beginning of a string only as a child
of the ``parent" string 1122..., so the child must begin with
2122... or 21102... (depending on length of first run of 2's in the
parent).  If the latter, then the child splits as 2110.2..., so the
new run of 2's is part of a new particle (the ``down quark", actually)
and so needn't concern us.  If the former, then the parent must
have been 112211... (again due to parsing rules) so the child begins
as 2122101 (which gives an immediate split) or as 212221 (which
means the grandchild gives an immediate split).  We can conclude
that we generate at most one new run of 2's at the beginning of our
string.

Referring back to  (*)  , we consider what happens over four steps.
Every interior run of 2's produces a split either to its left or
to its right, on step two, three, or four.  It's possible that one
of these splits picks up a new leading run of 2's by step 4, but by
that time the run of 2's to its right, if not the last run in the
string, has created its own split.  Thus, each new, split, substring
has at most two runs of 2's, except possibly the initial substring,
which might have three runs:  one created at the beginning of the
string, one from run G, and one from run I in the first case where
I has length two.  Let's consider the possibilities in this case:

\begin{tabular}[t]{cccccc}
 & G & H & I & J & ... \\
(possibly some 1's) & (2's)  &  (1's) & 22  & (1's) & ...
\end{tabular}

By writing down all possible initial configurations for the leading
1's, for G, and for H, (and there aren't too many of them), we can
verify that all cases will split into substrings with at most two
runs of 2's.  And the longest irreducible ancient string with only
two runs of 2's is 1112221112221110, of length $16$.

To sum up:  by Lemmas \ref{eli}, \ref{elii}, and \ref{eliii}, every
string eventually decays to an ancient string, which certainly
splits into irreducible ancient strings.  By the preceding, the
largest irreducible ancient string has length 16.  An irreducible
string will of course have at most one 0, and that, if it occurs,
in the final position.  And an irreducible ancient string that is
not itself common has no runs of length 4 or more.  Consequently,
every string eventually splits into a combination of common particles
and essential ancient strings.  So, if every essential ancient
string is eventually common, then {\em every} string is eventually
common, and that completes the Small World Theorem.
\end{proof}

\subsection* {The Cosmological Theorem}
\def\nota{\hat{\ }}
\def\alt{\,\bigm|\,}

The computer verification of the decay into common particles depends
on being able to recognize all points at which a ternary string can
be split, so our next goal is a proposition characterizing these
splitting points.  This is facilitated by the introduction of a
little notation.

We will use a variation of Conway's notation for describing patterns
a string might match.  Following Conway we use square brackets $[$
and $]$ to match the absolute beginning or end of a string; any
literal digit (0,1, or 2) matches itself; $\nota x$ matches anything
\emph{except} an $x$ (including possibly the absolute end of the
string); and a list separated by bars, as in $(x_1\alt x_2\alt\cdots\alt
x_n)$ matches any one of the $x_i$'s.

In this context, the symbol $\sim$ can be read as ``matches'', so
a statement like
$$X\sim[1(22\alt222)\nota2$$
means, ``$X$ is any string that begins with a 1, followed by either
22 or 222, followed by either nothing at all or by 
something that is not a 2.''  For example,
the strings 1221, 122201, or 122 satisfy this condition; 122221,
however, does not.

\begin{defn}
If $A_n\sim[(0|1)$ for all $n\ge0$ then we will call $A$ a {\it
forever-leading-2-free} string. 
\end{defn}

\begin{lem}\label{leading2free}
A non-empty, essential ancient string $A$ is forever-leading-2-free
if and only if $A$ has one of the following forms:
\begin{enumerate}
\item $A\sim{\tt[0}$
or
\item $A\sim{\tt \mathbf{\lbrack}1\left(0\alt11\alt2\nota2\alt222\alt\ ]\ \right)}$
\end{enumerate}
\end{lem}
 
\begin{proof} 
It is easy to verify that all of the cases listed will never produce
an initial 2 under the ternary decay operation.  For the converse,
if $A$ is essential ancient and forever-leading-2-free, then the
first run in $A$ clearly has to be one of $0$, $1$, or $111$.  ($A$
itself can't begin with a 2, and $A$ can't begin with a run of two
1's; all other possibilities are ruled out by the definition of
``essential ancient.'')

If the initial run in $A$ is a $0$ or $111$ it is clear that that
$A$ is forever-leading-2-free no matter what follows.  If the initial
run in $A$ is a $1$, then that $1$ is followed by either a $0$ or
a run of up to three $2$'s (or the end of the string).  Again, it
is easy to verify that any of these possibilities {\it except} a
$22$ makes $A$ forever-leading-2-free.
\end{proof}

\begin{rem}
We could introduce a similar definition of ``forever-initial-1-free''
strings, but it is very easy to check that the only forever-initial-1-free
strings are those of the form $A\sim[22Y$, where $Y$ is
forever-leading-2-free.
\end{rem}

This sets up the following characterization of splitting points:

\begin{splittingtheorem}
An essential ancient string $X=LR$ splits as $X\sim L.R$ if and
only if one of the following holds:

\begin{enumerate}
\item $L\sim 0]$ (it follows from maturity of $X$ that $R\sim [\nota0$).
\item $L\sim 1]$ and $R\sim[22Y$, where $Y$ is forever-leading-2-free.
\item $L\sim 2]$ and $R$ is forever-leading-2-free.
\end{enumerate}
\end{splittingtheorem}

\begin{proof}
This is really an immediate corollary of the preceding proposition
and the remark on initial-1-free strings.
\end{proof}

As long as we are able to identify all the splitting points, we can simply test all the necessary strings for their eventual decay into common particles.  Which strings are necessary?  By the previous section, it suffices to check the essential ancient strings.

\begin{alg}[for proving the Cosmological Theorem]
For each essential ancient string $A$,
\begin{enumerate}
\item Let $S:=A$. 
\item Apply the splitting theorem as many times as possible to
split $S$ into a list $T(1),\ldots,T(n)$ of irreducible strings.
\item If all the $T(i)$ strings are in the list of common strings,
then proceed to the next ancient string $A$, and return to step 1.
(If all ancient strings have been examined, then halt.)
\item Otherwise let $S:=S_1$ (that is, allow $S$ to decay by one
iteration) and return to step 2.
\end{enumerate}
If this algorithm halts, then the Cosmological Theorem is shown to
be correct -- and, happily, it does halt.  A Mathematica package
which implements all the ideas in this paper is available from the
second author's web page at \url{http://www.wlu.edu/~siehlerj/computing/}
along with a demonstration notebook that illustrates the package
commands and carries out the verification of the Cosmological
Theorem.

The investigator that wishes to check or reproduce our results independently
may find the following table (Table \ref{decaytimetable}) more helpful
than the simple affirmative answer ``It halts.''  The $n$-th row
of the table contains an enumeration of the essential ancient strings
of length $n$, broken down by the number of iterations it takes for
them to split into common particles.  For example, the ``8'' that
appears in the seventh row means there are eight essential ancient
strings of length 7 which require 5 iterations before splitting
entirely into common particles; those eight are:
\[ 1121122, 1122122, 1221121, 2112122, 2121121, 
2221121, 1121220, \hbox{ and } 2122120 \]
\end{alg}

\begin{table}[!hbt]
\small
\begin{tabular}{lr|rrrrrrrrrrr|r}
 && \multicolumn{10}{c}{Iterations required to split into common particles}&\\
&&  0&  1&  2&  3&  4&  5&  6&  7&  8&  9&  10& Total  \\
 \hline
&1&  1&  1&  0&  0&  0&  0&  0&  1&  0&  0&  0&3  \\ &2&  3&  1&  0&  0&  0&  1&  1&  0&  0&  0&  0&6  \\
&3&  5&  3&  0&  1&  0&  0&  0&  2&  1&  0&  0&12  \\
&4&  9&  5&  1&  0&  1&  2&  1&  2&  1&  0&  0&22  \\
&5&  10&  9&  1&  3&  2&  4&  2&  5&  3&  0&  1&40  \\
&6&  19&  16&  1&  6&  7&  5&  4&  7&  7&  2&  0&74  \\
String\!\!\!&7&  17&  33&  5&  18&  14&  8&  5&  18&  14&  3&  1&136  \\
Length\!\!\!&8&  32&  48&  11&  31&  27&  22&  15&  28&  28&  6&  2&250  \\
&9&  32&  92&  18&  70&  58&  30&  25&  61&  56&  14&  4&460  \\
&10&  57&  142&  40&  130&  109&  63&  50&  108&  110&  28&  9&846  \\
&11&  59&  243&  72&  258&  207&  116&  94&  217&  215&  57&  18&1556  \\
&12&  99&  386&  133&  473&  399&  221&  185&  394&  420&  115&  37&2862  \\
&13&  109&  639&  238&  898&  759&  392&  340&  767&  815&  232&  75&5264  \\
&14&  173&  1008&  440&  1627&  1418&  762&  653&  1417&  1575&  459&  150&9682  \\
&15&  199&  1638&  749&  3014&  2668&  1369&  1222&  2709&  3035&  906&  299&17808  \\
&16&  304&  2591&  1341&  5431&  4974&  2567&  2310&  5033&  5829&  1783&  591&32754  \\
\end{tabular}
\caption{Decay times for essential ancient strings}\label{decaytimetable}
\end{table}

We also point out that the number of essential ancient strings of
length $n\ge2$ with no final zero can be computed by the function
\[f(n)= 2\sum _{k=\left\lceil n/4\right\rceil }^n \left( \binom{n-1}{k-1}-\sum _{r=1}^{\left\lfloor n/4\right\rfloor } (-1)^{r+1} \binom{k}{r} \binom{n-3r-1}{k-1}\right),\]
which is left as a small exercise in basic counting principles; alternatively, the following recursion formula:
$$f(n)=f(n-1)+f(n-2)+f(n-3)$$
is easily seen to hold for $n\ge4$.
The total number of strings that should appear in the $n$-th row of the
table (with or without a final zero) is thus given by
$$f(n-1)+f(n),$$
a further check to ensure that all necessary strings have been
represented in the table.  This tabulation finishes the verification
of the Cosmological Theorem.

\subsection*{Growth and Frequency}

The following chart shows how each of our 24 particles decays into
other particles.  We note right away that if one starts with a
collection of fermions, then only fermions will be produced. A few
moments spent with the chart will also verify that once a single
fermion appears in a string, then eventually all eight will appear
simultaneously.  However, this is not the case with the bosons, and
it is quite possible for a string's descendents to contain only a
few bosons at any given time. The first five bosons in the loop at
the top will, of course, cycle through indefinitely, but most of
the bosons in the long vertical chain will appear only momentarily.
One can see now why we named them after the hypothetical supersymmetric
particles!

Recall that in Conway's world, all 92 of his elements eventually
appear in every sufficiently old (non-trivial) string.  In the base
2 world, only eight or nine of Silke's eleven common strings would
appear. The base 3 case is the best of both worlds. It's possible
to have an initial string (such as E or M or U) which decays into
collections of just the eight fermions, and it's also possible to
find an initial string (such as Nm.Nt.E.Ph.Ne.E.Gl.Ne.E.Wb.Ne.E.Zb)
which eventually decays into strings containing all 24 particles!

\setlength{\unitlength}{0.9em}
\noindent\begin{picture}(43,30)

\put(8,21){\makebox(0,1){\underbar{\bf Fermions}}}
\put(23,24){\makebox(0,1){\underbar{\bf Bosons}}}
\put(34,14){\makebox(0,1){\underbar{\bf Neutrinos}}}

\put(34,12){\makebox(0,1){\bf 22}}   
\put(34,12){\vector(0,-1){1}}
\put(34,10){\makebox(0,1){22}}
\qbezier[8](35,11)(36,11.5)(35,12)\put(35,12){\vector(-2,1){.1}}     

\put(31,6){\makebox(0,1){\bf 11110}}
\put(31,6){\vector(0,-1){1}}
\put(31,4){\makebox(0,1){11110}}
\qbezier[8](29,5)(28,5.5)(29,6)\put(29,6){\vector(2,1){.1}}     

\put(37,6){\makebox(0,1){\bf 11112}}
\put(37,6){\vector(0,-1){1}}
\put(37,4){\makebox(0,1){11112}}
\qbezier[8](39,5)(40,5.5)(39,6)\put(39,6){\vector(-2,1){.1}}     

\multiput(20,14)(0,-2){3}{\vector(0,-1){1}}

\put(8,18){\vector(0,-1){1}}
\put(8,18){\makebox(0,1){\bf 10}}
\put(8,16){\makebox(0,1){\bf 1110}}\put(8,16){\vector(1,-1){1}}\put(8,16){\vector(-1,-1){1}}
\put(8,14){\makebox(0,1)[r]{10\ }}\put(8,14){\makebox(0,1){\bf .}}\put(8,14){\makebox(0,1)[l]{\bf\ 110}}
\multiput(9,14)(0,-2){3}{\vector(0,-1){1}}
\put(9,12){\makebox(0,1){\bf 2110}}
\put(9,10){\makebox(0,1){\bf 122110}}
\put(9,8){\makebox(0,1){\bf 11222110}}\put(9,8){\vector(1,-1){1}}\put(9,8){\vector(-1,-1){1}}
\put(9,6){\makebox(0,1)[r]{2110\ }}\put(9,6){\makebox(0,1){\bf .}}\put(9,6){\makebox(0,1)[l]{\bf\ 22110}}
\put(10,6){\vector(0,-1){1}}
\put(10,4){\makebox(0,1){\bf 222110}}
\put(10,4){\vector(1,-1){1}}\put(10,4){\vector(-1,-1){1}}
\put(10,2){\makebox(0,1)[r]{10\ }}\put(10,2){\makebox(0,1){\bf .}}\put(10,2){\makebox(0,1)[l]{\ 22110}}
\qbezier[17](13,3)(14.5,4.5)(13,6)\put(13,6){\vector(-1,1){.1}}    

\qbezier[18](6,15)(5,16)(7,18)\put(7,18){\vector(1,1){.1}}   
\qbezier[29](6,7)(4,9)(7,12)\put(7,12){\vector(1,1){.1}}     
\qbezier[86](8,3)(0,11)(7,18)								 

\put(19,18){\vector(0,-1){1}}
\put(19,18){\makebox(0,1){\bf 12}}
\put(19,16){\makebox(0,1){\bf 1112}}\put(19,16){\vector(1,-1){1}}\put(19,16){\vector(-1,-1){1}}
\put(19,14){\makebox(0,1)[r]{10\ }}\put(19,14){\makebox(0,1){\bf .}}\put(19,14){\makebox(0,1)[l]{\bf\ 112}}
\multiput(20,14)(0,-2){3}{\vector(0,-1){1}}
\put(20,12){\makebox(0,1){\bf 2112}}
\put(20,10){\makebox(0,1){\bf 122112}}
\put(20,8){\makebox(0,1){\bf 11222112}}\put(20,8){\vector(1,-1){1}}\put(20,8){\vector(-1,-1){1}}
\put(20,6){\makebox(0,1)[r]{2110\ }}\put(20,6){\makebox(0,1){\bf .}}\put(20,6){\makebox(0,1)[l]{\bf\ 22112}}
\put(21,6){\vector(0,-1){1}}
\put(21,4){\makebox(0,1){\bf 222112}}
\put(21,4){\vector(1,-1){1}}\put(21,4){\vector(-1,-1){1}}
\put(21,2){\makebox(0,1)[r]{10\ }}\put(21,2){\makebox(0,1){\bf .}}\put(21,2){\makebox(0,1)[l]{\ 22112}}
\qbezier[17](24,3)(25.5,4.5)(24,6)\put(24,6){\vector(-1,1){.1}}						

\qbezier[33](19,2)(17,0)(13,0)														
\qbezier[33](13,0)(9,0)(7,2)														
\qbezier[86](7,2)(-1,10)(6,17)														
\qbezier[43](17,15)(13,16.5))(9,18)\put(9,18){\vector(-3,1){.1}}					
\qbezier[39](17,7)(14,9.5)(11,12)\put(11,12){\vector(-1,1){.1}}						
\qbezier[65](21.5,22.5)(16,22.5)(9,17)\put(9,17){\vector(-2,-1){.1}}				

\put(21,20){\makebox(0,1)[r]{\bf 2\ }}\put(21,20){\makebox(0,1){\bf .}}\put(21,20){\makebox(0,1)[l]{\bf\ 12221}}
\put(20,20){\vector(-1,-1){1}}

\put(24,19){\vector(-1,1){1}}
\put(23,21){\vector(1,1){1}}
\put(26,22){\vector(1,-1){1}}
\put(27,20){\vector(-1,-1){1}}

\put(27,20){\bf 1221}
\put(25,18){\makebox(0,1){\bf 112211}}
\put(24,22){\makebox(0,1)[r]{1110\ }}\put(24,22){\makebox(0,1){\bf .}}\put(24,22){\makebox(0,1)[l]{\bf\ 211}}
\end{picture}

\begin{proof}[Proof of Growth Theorem]
It should be clear from the chart that, for a string $A$ which
splits into common particles, (i) the number of neutrinos in the
iterates $A_n$ remains constant; (ii) the number of bosons in $A_n$
is bounded by a linear function of $n$; and (iii) the number of
fermions in $A_n$ grows exponentially, as long as there is at least
one fermion or boson present in $A$ itself.  Consequently, as $n$
gets large, the fraction of the string made up of bosons and/or
neutrinos goes to zero, and we can restrict our attention to the
rate at which fermions are produced.

The fermion information can be encoded in a square ``transition
matrix'' $\Lambda$, with one column and row for each fermion:

\[ \Lambda = \bordermatrix {
 & E& M& D& B& U& S& T& C \cr
E& 0 & 1 & 0 & 0 & 0 & 0 & 1 & 0 \cr
M& 1 & 0 & 0 & 0 & 0 & 0 & 0 & 0 \cr
D& 0 & 0 & 0 & 0 & 1 & 0 & 0 & 1 \cr
B& 0 & 0 & 0 & 0 & 0 & 0 & 1 & 1 \cr
U& 0 & 1 & 0 & 0 & 0 & 0 & 0 & 0 \cr
S& 0 & 0 & 1 & 0 & 0 & 0 & 0 & 0 \cr
T& 0 & 0 & 0 & 1 & 0 & 0 & 0 & 0 \cr
C& 0 & 0 & 0 & 0 & 0 & 1 & 0 & 0 \cr
}
\]

The entries in the $C$ column (standing for the Charm fermion, $11222110$),
 show that this particle produces one $D$ and one $B$ particle 
 ($2110$ and $22110$) upon decay, as is easily verified.  Now,
one can also check that the 14th power of our transition matrix
$\Lambda$ has all entries positive.  That means that any string $A$
containing even one fermion will, after 14 iterations, contain {\it
all} the fermions.  Moreover, the Perron-Frobenius theorem \cite{BR}
guarantees that $\Lambda$ has a dominant eigenvalue $\lambda$ which
is a positive real with largest magnitude among all the eigenvalues
of $\Lambda$,  and the growth rate of our strings under iteration
is asymptotic to $\lambda^n$.  A plot of the eigenvalues shows our
$\lambda$,
\begin{figure}[!hbt]
\begin{center}
    \includegraphics{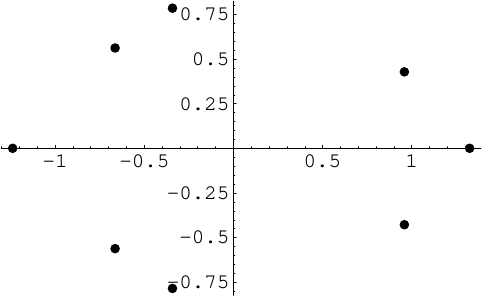}
    \caption{Eigenvalues of the fermion transition matrix}\label{fig-eigenvalue}
\end{center}
\end{figure}
and it is easy to compute (find and factor the characteristic
polynomial of $\Lambda$) that $\lambda$ is the unique positive
root of $x^3-x-1$, approximately $1.3247$.
\end{proof}

\begin{proof}[Proof of Frequency Theorem]
This is also derived directly from the transition matrix.  Computing
the limiting frequency with which a particular fermion occurs means
computing its row total in the powers $\Lambda^n$ of the transition
matrix, divided by the total of {\em all} entries in $\Lambda^n$, and
taking the limit as $n$ gets large.  Obtaining exact expressions
seems unlikely, but it is easy enough to produce the following
numerical approximations simply by considering ``large enough''
powers of $\Lambda$:

\begin{tabular}{l|l|l|l|l|l|l|l}
 E& M& D& B& U& S& T& C \\
 \hline
 18.50\% & 13.97\% & 13.97\% & 13.97\% & 10.54\% & 10.54\% & 10.54\% & 7.96\% \\
\end{tabular}

To be clear, this says that {\it any} string which contains
or produces even one non-neutrino particle will eventually be made
up of the eight fermions, occuring in the relative frequencies
above, together with (possibly) a vanishing amount of bosons and
neutrinos.  So we see that Charm, while not entirely absent, is
indeed the rarest of assets.
\end{proof}

\subsection*{For future research.} 
Let's first review what we have learned about the different
versions of the Look and Say sequences. The ``Range of $k$" appearing
in the chart refers to the function $k(S)$, which 
gives the number of different common particles that
appear in all sufficiently-old descendents of an initial string $S$. 
So, for example, $k(22) = 1$ in both Conway's world and in the base 3 world,
while for every other non-trivial string $S$, then in Conway's world $k(S) = 92$.

\V

\begin{tabular}{r|ccc}
 & Conway's base $\infty$ & base 3 & base 2 \\
 \hline
Elements & 92 & 24 & 11\\[1.2ex]
Fixed strings & {\bf 22} & {\bf 22}, {\bf 11112}, {\bf 11110} & {\bf 111} \\[1.2ex]
Range of $k$ & 1, 92 & 1, 2, 8-24 ? & 1,8,9 \\[1.2ex]
Growth rate &  $\lambda = 1.3035\dots$ & $\lambda = 1.3247\dots$ & $\lambda = 1.4655\dots$ \\[1.2ex]
Min. poly for $\lambda$ & (degree 71) & $x^3 - x - 1$ & $x^3 - x^2 -1$
\end{tabular}

\V

The reader will notice that we have not yet fully determined which
$k$ values between $1$ and $24$ are possible in the base 3 world.
We gave a (fairly long) string earlier that eventually decays into
all 24 particles, and we have examples with $k$ values of 1, 2, 8,
and 9, but we have not investigated the matter any further. Perhaps
there are some strings with no well-defined $k$ value?

Changing the nature of the problem (and based on a suggestion from
David Park), Conway has investigated a ``Look and Say" sequence  of
Roman numbers: I, then I I, followed by II I, then III I, and IV
I, and so on.  He states in  \cite{C2} that the $\lambda$ associated
with this sequence is not of particularly high degree, but gives
no further details.  We wonder if this sequence behaves more like
the base 10 or the trinary system. We also wonder if there are some
other ancient numbering systems (such as the Mayan or Babylonian
numbers) that might be worth looking into.

One could also study an alphabetic version of these sequences:

\begin{center}
One.\\
One e, one n, one o.
\end{center}

If we think of this as a ``Counting Sequence", the next line would be:

\begin{center}
four e's, four n's, four o's
\end{center}

On the other hand, if we think of this as a ``Look and Say" sequence,
the next line would be
\begin{center}
One o, one n, two e's, one o, one n, one e, one n, one o, one n, one e, one o.
\end{center}
It's interesting to note that the ``Counting Sequence" method does
admit fixed strings; for example:
\begin{quote}
One~a, one~b, one~d, twenty-four~e's, six~f's, four~g's, seven~h's, eight~i's, one~l, thirteen~n's, thirteen~o's, one~q, eight~r's, eighteen~s's, fourteen~t's, five~u's, three~v's, four~w's, two~x's, two~y's, one z.
\end{quote}
Of course, the structure of these alphabetic counting sequences is language-dependent.  An example in Hungarian (well, why not?) verifies that other languages do admit fixed strings:
\begin{quote}
Egy~a, \"ot~\'a, egy~b, n\'egy~c, egy~d, tizenk\'et~e, nyolc~\'e, tizenh\'arom gy, h\'et~h, n\'egy~i, egy~\'i, egy~j, h\'arom~k, n\'egy~l, \"ot~m, nyolc~n, h\'arom ny, h\'et~o, n\'egy~\"o, egy~\H o, \"ot~r, egy sz, kilenc~t, egy~u, h\'arom z.
\end{quote}
Note that {\it gy} and {\it ny} are digraphs which are treated as a single letter in the Hungarian alphabet; also, differently-accented vowels such as o, \"o, and \H o, are counted as distinct from one another.  Is the existence of these strings language-independent, or is there some strange language which would never admit such a self-documenting string?

As for the ``Look and Say" method of interpreting these alphabetical
strings, it would be extremely difficult to determine the common elements and 
the growth rate; we wouldn't want to attack this without a computer! It
would be interesting to determine if different languages gave a different
growth rate $\lambda$.

\subsection*{Appendix.} 

In Conway's original article, he named his 92 common strings after
the elements of the periodic table. In our base-3 world, we have 
chosen to name our
24 elements after the sub-atomic particles.
The following brief explanation is not necessary for
understanding the mathematics in this paper, but it might help to
illustrate why we selected these particular names for our strings.
Recall that in quantum field theory, the {\bf elementary particles}
are objects (like quarks, electrons, photons, and so on) which are
believed to  have no internal structure and thus can not be divided
into smaller objects (just like our 24 strings).
In physics, these elementary particles can combine into the more
familiar protons and neutrons (among other things), which then form
atoms and molecules, just as our 24 strings can combine to form
longer strings.

As seen in Table \ref{base3}, we group our particles into {fermions},
{bosons}, and {neutrinos}. The first group of eight we call fermions,
and it consists
of two leptons (in physics, the familiar
electron and the muon) followed by  six quarks (in physics,
given the names Up, Down, Strange, Charm, Bottom, and Top).
The next column of thirteen  particles we denote as {bosons};
the first five are the somewhat familiar photon, gluon, W and Z
bosons, and Higgs boson, while the next eight are denoted as the
supersymmetric partners of the leptons and quarks in the first
column. Note that these supersymmetric bosons are almost identical
to the fermions to their left, differing only in the last digit.
Their names reflect this symmetry; 222112 is quite similar to the
``top quark" 222110 (T), and so we named it after the ``stop squark"
(St) subatomic particle.  In physics these supersymmetric particles
are hypothetical, and likewise (in our Frequency Theorem above)
we have shown that the eponymous
strings appear in a vanishingly small proportion,  if at all,
under iteration of the Look and Say operation.
Finally, our last column
in Table \ref{base3}
identifies three ``boring" strings as {neutrinos} (the electron
neutrino, the mu neutrino, and the tau neutrino). These three
strings are fixed under the Look-and-Say operation, and so we feel
they play a ``neutral'' role in the dynamics.

\V\V\V\V


%
%
%

\end{document}